\documentclass[a4,psfig,english,10pt,epsf,portrait]{article}
\usepackage{amsmath,amsfonts,latexsym,amscd,amssymb,theorem}
\usepackage{epsfig}
\usepackage{amsmath}
\def\R{\hbox{\bf R}}
\def\Z{\hbox{\bf Z}}

\def\o{\overline}

\def\N{\hbox{\bf N}}

\def\g{\gamma}

\def\a{\alpha}

\def\D{{\cal D}}

\def\l{{\lambda}}

\def\D{\Delta}
\def\sgn{\mbox{sgn}}

\def\t{\tau}

\def\<{\langle}
\def\>{\rangle}
\newcommand{\ba}{\begin{eqnarray}}
\newcommand{\ea}{\end{eqnarray}}

\newtheorem{thm}{Theorem}[section]

\newtheorem{theorem}[thm]{Theorem}

\newtheorem{lemma}[thm]{Lemma}
\newtheorem{proposition}[thm]{Proposition}
\newtheorem{corollary}[thm]{Corollary}
\newtheorem{rem}[thm]{Remark}
\newtheorem{example}[thm]{Example}

\newcommand{\eps}{\epsilon}
\numberwithin{equation}{section}

\renewcommand{\R}{{\mathbb R}}
\renewcommand{\Z}{{\mathbb Z}}
\renewcommand{\N}{{\mathbb N}}

\renewcommand{\g}{\gamma}

\begin{document}

\title{\bf On the rate of convergence in periodic homogenization of
  scalar first-order ordinary differential equations}
\author{
\normalsize\textsc{ H. Ibrahim $^{*,\dagger}$, R. Monneau
  \footnote{Universit\'{e} Paris-Est,
Ecole des Ponts, CERMICS, 6 et 8 avenue Blaise Pascal, Cit\'e
Descartes Champs-sur-Marne, 77455 Marne-la-Vall\'ee Cedex 2, France.
E-mails: ibrahim@cermics.enpc.fr, monneau@cermics.enpc.fr 
\newline \indent $\,\,{}^\dagger$CEREMADE, Universit\'{e}
Paris-Dauphine, Place De Lattre de Tassigny, 75775 Paris Cedex 16, France} }}
\vspace{20pt}

\maketitle


\centerline{\small{\bf{Abstract}}} \noindent{\small{In this paper, we
    study the rate of convergence in periodic homogenization of scalar
    ordinary differential equations. We provide a quantitative
    error estimate between the solutions of a first-order ordinary
    differential  equation with rapidly oscillating coefficients and the
    limiting homogenized solution. As an application of our
    result, we obtain an error estimate for the solution of some
    particular linear transport equations.}} \hfill\break

%
%
%
%
%
%

\section{Introduction}
\subsection{Homogenization of an ODE}
In this paper, we consider the solutions of the following first-order
ordinary differential equation:
\begin{equation}\label{main_eq}
\left\{
\begin{aligned}
&u^{\eps}_{t}=f\left(\frac{u^{\eps}}{\eps},\frac{t}{\eps},u^{\eps},t\right),\quad
t>0,\\
&u^{\eps}(0)=u_{0},
\end{aligned}
\right.
\end{equation}
where $\eps>0$, $u^{\eps}_{t}$ stands
for $\frac{du^{\eps}}{dt}$ or equivalently $\partial_{t}u^{\eps}$, and
$u_{0}$ is a real number. We are interested in the rate of convergence
of the solution $u^{\eps}$ to its limit in the framework of periodic
homogenization. We employ the following assumptions on
the function $f$:
\begin{itemize}
\item (A1) {\bf Regularity:} the function $f:\R^{4}\rightarrow \R$ is a
  bounded Lipschitz continuous function with $\alpha=Lip(f)$
its Lipschitz constant, and $\beta=\|f\|_{L^{\infty}(\R^{4})}$;

\item (A2) {\bf Periodicity:} for any $(v,\tau,u,t)\in \R^{4}$, we
have:
$$f(v+l,\tau+k,u,t)=f(v,\tau,u,t)\quad \mbox{for any}\quad  (l,k)\in \Z^{2};$$

\item (A3) {\bf Monotonicity:} for any $(v,\tau,t)\in \R^{3}$,
$$ \mbox{the function } u\longmapsto f(v,\tau,u,t) \mbox{ is non-increasing}.$$

\end{itemize}
Let us make short comments on these assumptions. Remark that assumption (A1) ensures
the existence and uniqueness of the solution $u^{\eps}$ of
(\ref{main_eq}) via the Cauchy-Lipschitz theorem. Moreover, the assumed
boundedness of $f$ is not a restrictive condition while we work on any
finite time interval $[0,T]$. The monotonicity assumption (A3) may seem
unnecessary at a first glance, but will be indeed useful to guarantee
the uniqueness of the solution to the homogenized equation (see
Proposition \ref{ganaza}). Moreover, assumption (A3) will play a crucial
role to establish the rate of convergence of $u^{\eps}$ to its limit
$u^{0}$ (see for instance Section \ref{sec4}).

In order to define the homogenized equation, we will use the following
proposition:
\begin{proposition}\textbf{(Definition and properties of the effective
    slope $\o{f}$).}\label{AM7Z}
Fix $(u,t)\in \R^{2}$. Then there exists $\l\in \R$ such that for any
initial data $u_{0}\in \R$, the solution $v\in C^{1}([0,\infty);\R)$ of
the following ordinary
differential equation:
$$
\left\{
\begin{aligned}
& v_{\tau}=f(v, \t, u, t),\quad \t>0,\\
& v(0)=u_{0},
\end{aligned}
\right.
$$
satisfies
\begin{equation}\label{Ymala1}
\frac{v(\t)}{\t}\rightarrow \l\quad \mbox{as}\quad \t\rightarrow \infty.
\end{equation}
Let us set the effective slope:
\begin{equation}\label{Ymala2}
\overline{f}(u,t)=\l.
\end{equation}
Then the following holds:
\begin{equation}\label{Ymala3}
\left\{
\begin{aligned}
& \o{f}: \R^{2}\rightarrow \R\mbox{ is continuous.}\\
& \mbox{For any $t\geq 0$, the map $u \mapsto \o{f}(u,t)$ is
  non-increasing.}
\end{aligned}
\right.
\end{equation}
\end{proposition}
Let us mention that, for some specific functions $f$, explicit formulas
for $\o{f}$ can be obtained (see for instance \cite{Picc78}, and the
examples below).
\begin{example}\label{example1}
For $f(v,\t,u,t)=-u+\cos(2\pi v)$, we have
$\o{f}(u,t)=\left(\int^{1}_{0}\frac{dv}{-u + \cos(2\pi v)}
\right)^{-1}$ for $u>1$,
and
$$\o{f}(u,t)\sim -c\sqrt{u-1} \quad\mbox{as}\quad u \rightarrow
1^{+},$$
for some constant $c>0$.
\end{example}
\begin{example}\label{example2}
For $f(v,\t,u,t)=-u+|\sin(2\pi v)|$, we have, for some constant $c>0$:
\begin{equation}\label{ineq_example}
\o{f}(u,t)\sim \frac{c}{|\log |u||}\quad \mbox{as}\quad
u\rightarrow 0^{-}.
\end{equation}
\end{example}
Example \ref{example1} shows in particular that even for analytic $f$,
the function $\o{f}$ could be non-Lipschitz. Example~\ref{example2}
shows a case where $\o{f}$ is not H\"{o}lder continuous when $f$ is
Lipschitz continuous. The proof of Example~\ref{example2} will be given
in the Appendix. At this stage, we can write the homogenized equation associated to
equation (\ref{main_eq}) as follows:
\begin{equation}\label{homog_eq}
\left\{
\begin{aligned}
& u^{0}_{t}=\overline{f}(u^{0},t),\quad t>0,\\
& u^{0}(0)=u_{0}.
\end{aligned}
\right.
\end{equation}
Even if $\o{f}$ may not be Lipschitz continuous in $u$, we can show the
existence and uniqueness of the solution of (\ref{homog_eq}), taking
advantage of the monotonicity of $\o{f}(u^{0},t)$ in $u^{0}$. Indeed, we have:
\begin{proposition}\label{ganaza}\textbf{(Existence and uniqueness).}
Under assumption (\ref{Ymala3}) on $\o{f}$, there exists a unique
solution $u^{0}\in C^{1}([0,\infty);\R)$ of (\ref{homog_eq}).
\end{proposition}
It is worth noticing that assumption (\ref{Ymala3}) satisfied by
$\o{f}$ in the homogenized equation is our motivation to make
assumption (A3) on $f$. We can now state our main result:
\begin{theorem}\textbf{(Error estimate for ODEs).}\label{theo2}
Under assumptions (A1)-(A2)-(A3), if $u^{\eps}$ is the solution of
(\ref{main_eq}), and $u^{0}$ is the solution of the homogenized
equation (\ref{homog_eq}), then for every $C>0$, $\eps>0$, and every
$T\geq C \eps |\log\eps|$, we have the following estimate:
\begin{equation}\label{error}
\|u^{\eps}-u^{0}\|_{L^{\infty}(0,T)}\leq \frac{c T}{|\log\eps|},
\end{equation}
where $c>0$ is a positive constant only depending on $C$ and on $\a$,
$\beta$ defined in assumption (A1).
\end{theorem}
Such a result for a general monotone system of ODEs seems to be
completely open. The above estimate in $\frac{1}{|\log \eps|}$ is in
fact related to the behavior of $\o{f}$ in Example \ref{example2}, which
is the worst possible regularity of $\o{f}$. Moreover, it is possible to
show that under the condition $T\geq C\eps |\log \eps|$, inequality (\ref{error})
is sharp, see the following example whose proof will be given in the
Appendix: 
\begin{example}\label{example3}
Let $f(v,\t,u,t)=g(v+\t)-1$ with a 1-periodic function $g$ satisfying
\begin{equation}\label{tresb}
g(w)=|w-{1}/{2}|\quad \mbox{for}\quad w\in [0,1].
\end{equation}
In this case $\o{f}(u,t)=-1$. Let us choose the initial data $u_{0}=0$. Then
for any $\delta >0$, we have the following estimate between the solution
$u^{\eps}$ to (\ref{main_eq}) and $u^{0}$ to (\ref{homog_eq}):
$$u^{\eps}(t)-u^{0}(t)\sim \frac{t}{2 \delta |\log \eps|}\quad
\mbox{for}\quad  t=\delta\eps |\log \eps|.$$ 
\end{example}
\begin{rem}
It is worth mentioning that assumption (A1) could be replaced by the
weaker assumption:
\begin{itemize}
\item {\rm{(A1){'}}} {\bf Regularity:} the function $f: \R^{4}\rightarrow
\R$ is a bounded continuous function such that for every $\t\in \R$,
the function $f(.,\t,.,.)$ is Lipschitz continuous.
\end{itemize}
\end{rem}
\subsection{Application to the homogenization of linear transport
  equations}\label{sub1.2} 
For $x=(x_{1},x_{2})\in \R^{2}$, let us consider a vector field
$a^{\eps}=(a^{\eps}_{1},a^{\eps}_{2})$ defined as follows
\begin{equation}\label{att_sign}
\left\{
\begin{aligned}
& a_{1}^{\eps}(x_{1},x_{2})=- f\left(\frac{x_{1}}{\eps},
  \frac{x_{2}}{\eps}, x_{1}, x_{2}\right)\\
& a^{\eps}_{2}(x_{1},x_{2})=1,
\end{aligned}
\right.
\end{equation}
with a function $f$ satisfying (A1)-(A2)-(A3). We consider the viscosity
solution $V^{\eps}(t,x)$ of the following linear transport equation:
\begin{equation}\label{transport_eqn}
\left\{
\begin{aligned}
& V^{\eps}_{t}+a^{\eps}\cdot \nabla V^{\eps}=0\quad &\mbox{on}&\quad
(0,\infty)\times \R^{2}\\
& V^{\eps}(0,x)=V_{0}(x)\quad &\mbox{on}&\quad \R^{2},
\end{aligned}
\right.
\end{equation}
where $V_{0}:\R^{2}\rightarrow \R$ is a Lipschitz continuous
function. The existence and uniqueness of a viscosity solution
$V^{\eps}$ of (\ref{transport_eqn}) is ensured since $a^{\eps}\in
W^{1,\infty}(\R^{2})$ and $V_{0}$ is Lipschitz continuous (see
for instance \cite{Barles}). The expected homogenized equation
associated to
(\ref{transport_eqn}) is:
\begin{equation}\label{transport_eqn0}
\left\{
\begin{aligned}
& V^{0}_{t}+\overline{a}\cdot \nabla V^{0}=0\quad &\mbox{on}&\quad
(0,\infty)\times \R^{2}\\
& V^{0}(0,x)=V_{0}(x)\quad &\mbox{on}&\quad \R^{2},
\end{aligned}
\right.
\end{equation}
with the vector field $\o{a}=(\o{a}_{1},\o{a}_{2})$ defined as:
\begin{equation}\label{att_sign_bar}
\left\{
\begin{aligned}
& \o{a}_{1}(x_{1},x_{2})=- \o{f}(x_{1}, x_{2}),\quad
\mbox{with $\o{f}$ given by (\ref{Ymala2})}.\\
& \o{a}_{2}(x_{1},x_{2})=1.
\end{aligned}
\right.
\end{equation}
As a consequence of Theorem \ref{theo2}, we will show in Section
\ref{sec6} the following result:
\begin{theorem}\label{theo3}\textbf{(Error estimate for linear transport
    equations).}
Under the previous assumptions, there exists a Lipschitz continuous
function $V^{0}$ which is a viscosity solution of
(\ref{transport_eqn0}), such
that for any $C>0$, $\eps>0$, and $T\geq C\eps|\log \eps|$, the solution
$V^{\eps}$ of (\ref{transport_eqn}) satisfies:
\begin{equation}\label{1.12}
\|V^{\eps}-V^{0}\|_{L^{\infty}(\R^{2}\times (0,T))}\leq \frac{c'T}{|\log
  \eps|}
\end{equation}
where $c{'}=c Lip(V_{0})$, and $c>0$ is the constant given in
Theorem \ref{theo2}.
\end{theorem}
Choosing the initial condition $V_{0}(x)=x_{1}$, we can easily deduce
from Example \ref{example3} that inequality (\ref{1.12}) is also sharp
for $T\geq C\eps |\log \eps|$. Here, in this application, the vector
field $a^{\eps}$ is quite special. The interested reader could be
referred to \cite{HouXin} for some other examples of vector fields in
$2$D, where a homogenization result is presented without any rate of
convergence. In \cite{Tassa97}, the author
gives some non-explicit error estimates for linear transport equations
in the particular case of periodic vector field $a^{\eps}$. However,
these error estimates obtained in \cite{Tassa97} may depend strongly on the
irrationality of the rotation number $\omega_{0}$ associated to the vector
field $a^{\eps}$ (where $\omega_{0}$ is nothing else than $-\o{f}$ in our
application). On the contrary, estimate (\ref{1.12}) only depends on
some bounds of the data of the problem, and are completely uniform with
respect to the rotation number.

Remark that when $f(v,\t,u,t)$ is independent of $u$ and $t$, we have a
much better estimate:
\begin{theorem}\label{theo4}\textbf{(Better error estimate).}
Under the assumptions of Theorem \ref{theo3}, if $f(v,\t,u,t)$ is
independent of $u$ and $t$, then we have:
$$\|V^{\eps}-V^{0}\|_{L^{\infty}(\R^{2}\times (0,T))}\leq c{''}
\eps,\quad \forall T\geq 0,$$
where $c{''}=\xi Lip(V_{0})$, with $\xi$ (given in Proposition
\ref{ergo}) only depends on $\beta$ defined in assumption (A1).  
\end{theorem} 
The proof of Theorem \ref{theo4} will also be given in Section \ref{sec6}.
\subsection{Brief review of the literature}
The pioneering work (via the theory of viscosity solutions) to periodic
homogenization was established in \cite{LPV_UnP}. Starting from
\cite{LPV_UnP}, the homogenization theory for Hamilton-Jacobi equations
has received a considerable interest. There is a huge literature that we
cannot cite in details, but the interested reader can for instance see
\cite{Alvarez-Bardi03, Barles-Soug01, Barles-Soug00, Evans92,
  Lions-Souga05, Hom-MonneauI, Hom-MonneauII} and the references therein. Another aspect
concerning homogenization of SDEs (stochastic differential equations)
has also been studied by several authors (see for instance
\cite{Pardoux99, GaPa01, BePa07, SOUG99}). These problems are related to
our problem when the SDE reduces to an ODE.

To our knowledge, the question of estimating the rate of convergence in
homogenization of PDEs has not been widely tackled up elsewhere in the
literature. We can cite \cite{CCG08} for several error estimates
concerning the rate of convergence of the approximation scheme to the
effective Hamiltonian. We can also cite the
work in \cite{C.D.-I} about the rate of convergence in periodic
homogenization of first-order stationary Hamilton-Jacobi 
equations, where an error estimate in $\eps^{1/3}$ is obtained for
Hamilton-Jacobi equations with Lipschitz effective Hamiltonian.

For the problems of homogenization of ODEs, we refer the reader to
\cite{Picc78,Picc79}. We also refer the reader to \cite{AHZ94,
Dalibard, WeinanE92, Menon02, Peirone96, Petrini99, Tartar89} for problems on
homogenization of nonlinear first-order ODEs and/or the associated
linear transport equations. As mentioned above, we refer the reader
to \cite{Tassa97} for some other error estimates for linear
transport equations.
\subsection{Organization of the paper}
The paper is organized as follows. In Section \ref{sec2}, we present
the proof of an ergodicity result (Proposition \ref{ergo})
that defines $\o{f}=\l$. We also present the proofs of Propositions \ref{AM7Z}
and \ref{ganaza}. In Section \ref{sec3}, we give a result of stability of $\l$
under additive perturbation (Proposition \ref{prop_pert}). A basic error
estimate (Proposition \ref{bef_proof}) is presented in Section
\ref{sec4}. Section \ref{sec5} is devoted to show our main result of
estimating the rate of convergence (Theorem \ref{theo2}). In
Section \ref{sec6}, we give an application to the case of linear
transport equations (Theorems \ref{theo3} and \ref{theo4}). We end up
in Section \ref{sec7} with an Appendix where we give the proof of
Examples \ref{example2} and \ref{example3}.

%
%
%
%
%
%
%

\section{Ergodicity and preliminary facts}\label{sec2}
In this section we present the proof of Propositions \ref{AM7Z} and
\ref{ganaza}. We first start with the following ergodicity
result which is a particular case of \cite[Proposition 4.2]{FIM}. 
However, we give the proof in our particular case for the
sake of completeness. 
\begin{proposition}\textbf{(Ergodicity).}\label{ergo} Let
  $g(v,\tau):\R^{2}\rightarrow \R$ be a
function satisfying:
\begin{itemize}
\item (H1) {\bf Regularity:} $g$ is Lipschitz continuous with
  $\|g\|_{L^{\infty}(\R^{2})}\leq \beta$;
\item (H2) {\bf Periodicity:} $g(v+l,\t+k) = g(v, \t)$ for any $(l,k)\in
  \Z^{2}$, $(v,\t)\in \R^{2}$.
\end{itemize}
\noindent Let $v$ be the solution of the following equation 
\begin{equation}\label{ergo_eq}
\left\{
\begin{aligned}
& v_{\tau}=g(v,\tau),\quad \tau>0,\\
& v(0)=v_{0},
\end{aligned}
\right.
\end{equation}
then there exist a constant $\lambda\in \R$ (independent of the initial
data $v_{0}$) such that for every $\t, \t'\geq 0$, we have:
\begin{equation}\label{ergo_eq1}
|v(\t)-v(\t{'})-\lambda(\tau-\tau')|\leq \xi \quad\mbox{with}\quad
\xi=1+2\beta. 
\end{equation}
\end{proposition}
\begin{rem}
Under our assumptions, it is possible (see \cite[Theorem 1.5]{FIM}) to
show the existence of a hull function $h:\R\times \R \rightarrow \R$
satisfying:
$$
\left\{
\begin{aligned}
& h(\t+1,x)=h(\t,x)\\
& h(\t,x+1)=h(\t,x)+1\\
& h_{x}\geq 0,
\end{aligned}
\right.
$$ 
such that $U(\t,x)=h(\t, x+\l\t)$, with $\l$ given in Proposition
\ref{ergo}, is a viscosity solution of:
$$ 
\left\{
\begin{aligned}
& U_{\t}=g(U,\t)\\
& U(\t=0, x)=x.
\end{aligned}
\right.
$$
This function $h$ may be discontinuous, but its existence suggests that
\begin{equation}\label{doble_etoi}
v(\t)=h(\t, \l\t)
\end{equation}
is formally a classical of (\ref{ergo_eq}). The expression
(\ref{doble_etoi}) allows to understand estimate (\ref{ergo_eq1}) and
also suggests that $v(\t)-\l\t$ could be quasi-periodic in some
cases. We emphasize the fact that in the proof of Proposition
\ref{ergo}, we do not use any notion of hull function, but propose a
completely independent proof.
\end{rem}
\textbf{Proof of Proposition \ref{ergo}.} For any $T>0$, define the two quantities:
$$\l^{+}(T)=\sup_{\tau\ge 0}\frac{v(\tau+T)-v(\tau)}{T} \quad\mbox{and}\quad
\l^{-}(T)=\inf_{\tau\ge 0}\frac{v(\tau+T)-v(\tau)}{T}.$$ These
quantities are finite since $v_{\t}=g$ is bounded. The proof is divided
into three steps.\\

\noindent {\bf Step 1: Estimate of $|\l^{+}-\l^{-}|$}.\\

\noindent Let $\delta>0$ be an arbitrary constant. From the
definition of $\l^{\pm}(T)$, there always exists $\t^{\pm}$ such that
\begin{equation}\label{deltaest}
\left|\l^{\pm}(T)-\frac{v(\t^{\pm}+T)-v(\t^{\pm})}{T} \right|\leq
\delta.
\end{equation}
Denote by $\lfloor t \rfloor$ and $\lceil t \rceil$ as the floor
and the ceil integer parts of the real number $t$ respectively. We
consider
\begin{equation}\label{shu0}
k=\lfloor \t^{-}-\tau^{+} \rfloor,\quad
\widetilde{\t}^{-}=\t^{-}-k\quad \mbox{and}\quad l=\lceil
v(\widetilde{\t}^{-})-v(\t^{-}) \rceil,
\end{equation}
and consider $w(\tau)=v(\tau+k)+l$. Using (\ref{ergo_eq}) and (H2), we
check that $w$ is a solution of the following equation
\begin{equation}\label{eqn_w}
\left\{
\begin{aligned}
& w_{\t}=g(w,\t),\quad \tau >0\\
& w(\widetilde{\t}^{-})=v(\t^{-})+l.
\end{aligned}
\right.
\end{equation}
We remark, from (\ref{shu0}), that
\begin{equation}\label{1}
\t^{+}\leq \widetilde{\t}^{-}< \t^{+}+1,
\end{equation}
and
\begin{equation}\label{2}
v(\widetilde{\t}^{-})\leq w(\widetilde{\t}^{-})<
v(\widetilde{\t}^{-})+1.
\end{equation}
From  (\ref{ergo_eq}), (\ref{eqn_w}) and (\ref{2}), the comparison
principle for ODEs gives in particular 
\begin{equation}\label{chiWs}
v(\tau)\leq w(\tau)\quad \mbox{for all}\quad \tau\geq
\widetilde{\t}^{-}.
\end{equation}
If we suppose $T\geq 1$, we obtain from (\ref{1}) that
$\t^{+}+T\geq \widetilde{\t}^{-}$ and hence (using (\ref{chiWs})):
\begin{equation}\label{3}
v(\t^{+}+T)\leq w(\t^{+}+T).
\end{equation}
Direct computations give
\begin{eqnarray*}
(v(\t^{-}+T)-v(\t^{-}))-(v(\t^{+}+T)-v(\t^{+}))
&=&(w(\widetilde{\t}^{-}+T)-w(\t^{+}+T))+(v(\widetilde{\t}^{-})-w(\widetilde{\t}^{-}))\\
& &\hspace{-0.3cm}+(w(\t^{+}+T)-v(\t^{+}+T))+(v(\t^{+})-v(\widetilde{\t}^{-})),
\end{eqnarray*}
where, from (\ref{2}) and (\ref{3}), we deduce that $
(v(\t^{-}+T)-v(\t^{-}))-(v(\t^{+}+T)-v(\t^{+}))\geq
-1-2\|g\|_{\infty}$. This inequality, together with (\ref{deltaest})
show that $ 0\leq \l^{+}(T)-\l^{-}(T)\leq
\frac{1+2\|g\|_{\infty}}{T}+2\delta, $ and since this is true for
any $\delta>0$, we obtain
$$|\l^{+}(T)-\l^{-}(T)|\leq \frac{1+2\|g\|_{\infty}}{T}.$$
However, in the case where $T\leq 1$, we always have
$\frac{|v(\t+T)-v(\t)|}{T}\leq \|g\|_{\infty} \leq \frac{\|g\|_{\infty}}{T}$
and therefore
$$0\leq \l^{+}(T)-\l^{-}(T)\leq \frac{2\|g\|_{\infty}}{T},$$
then
\begin{equation}\label{fin_est_lam}
|\l^{+}(T)-\l^{-}(T)|\leq \frac{\xi}{T}\quad 
\mbox{for every } T>0.
\end{equation}

\noindent {\bf Step 2: Existence of the limit of
$\l^{\pm}$(T) as $T\rightarrow \infty$}.\\

\noindent First, if we compute $\l^{+}(PT)$ for
$P\in\N\setminus\{0\}$, we get:
$$\l^{+}(PT)=\sup_{\tau>0}\frac{1}{P}\left[\sum_{i=1}^{P}\frac{v(\t+iT)-v(\t+(i-1)T)}{T}\right]\leq
\l^{+}(T).$$ Similarly, we get $\l^{-}(PT)\geq \l^{-}(T)$. Consider
$T_{1}$, $T_{2}>0$  such that
$T_{1}P=T_{2}Q$ for some $P, Q\in \N\setminus\{0\}$. Using this and
(\ref{fin_est_lam}), we have 
$$\l^{+}(T_{2}) \geq\l^{+}(T_{2}Q) =\l^{+}(T_{1}P)\geq
\l^{-}(T_{1}P) \geq \l^{-}(T_{1}) \geq
\l^{+}(T_{1})-\frac{\xi}{T_{1}},$$
similarly we have
$\l^{+}(T_{2})-\l^{+}(T_{1})\leq
\frac{\xi}{T_{2}},$ then
\begin{equation}\label{lambda+}
|\l^{+}(T_{1})-\l^{+}(T_{2})|\leq
\max\left(\frac{\xi}{T_{1}},
\frac{\xi}{T_{2}}\right).
\end{equation}
By the same arguments as above we can get
\begin{equation}\label{lambda-}
|\l^{-}(T_{1})-\l^{-}(T_{2})|\leq
\max\left(\frac{\xi}{T_{1}},
\frac{\xi}{T_{2}}\right).
\end{equation}
Recall that (\ref{lambda+}) and (\ref{lambda-}) are true when
$T_{1}/T_{2}$ is rational. By an approximation argument, joint with the
continuity of $\l^{\pm}$, it is easy to see that this is still true when
$T_{1}/T_{2}$ is any positive real number. Moreover, the identities
(\ref{lambda+}) and (\ref{lambda-}) give that the 
sequence $(\l^{\pm}(T))_{T}$ is a Cauchy sequence as $T\rightarrow
\infty$, and hence it has a limit:
\begin{equation}\label{l+}
\lim_{T\rightarrow \infty}\l^{\pm}(T)=\l,
\end{equation}
which is the same limit because of (\ref{fin_est_lam}).
From (\ref{l+}), (\ref{lambda+}) and (\ref{lambda-}), inequality
(\ref{ergo_eq1}) directly follows.\\

\noindent {\bf Step 3: Independence of $v_{0}$}.\\

\noindent The fact that $\l$ is independent of $v_{0}$ follows directly
from the comparison principle and inequality (\ref{ergo_eq1}). $\hfill{\Box}$\\

We can now present the proof of Propositions \ref{AM7Z} and \ref{ganaza}.\\

\noindent {\bf Proof of Proposition \ref{AM7Z}.} Using inequality
(\ref{ergo_eq1}) of Proposition \ref{ergo}, we can easily see that
(\ref{Ymala1}) directly follows. It remains to show (\ref{Ymala3}). We
argue in two steps.\\

\noindent {\bf Step 1: Monotonicity of $\o{f}$}.\\ 

\noindent Let $u\leq \tilde{u}$. Call $\l=\o{f}(u,t)$ and
$\tilde{\l}=\o{f}(\tilde{u},t)$. Let $v$ and $\tilde{v}$ be the solutions of:
$$
\left\{
\begin{aligned}
& v_{\tau}=f(v,\tau,u,t),\quad \tau>0,\\
& v(0)=u_{0},
\end{aligned}
\right.
$$
and
$$
\left\{
\begin{aligned}
& \tilde{v}_{\tau}=f(\tilde{v},\tau,\tilde{u},t),\quad \tau>0,\\
& \tilde{v}(0)=u_{0},
\end{aligned}
\right.
$$
respectively. Assume without loss of generality that $u_{0}=0$.
Using (A3), we deduce that
$f(\tilde{v},\tau,\tilde{u},t)\leq f(\tilde{v},\tau,u,t)$. Hence, the comparison principle gives:
\begin{equation}\label{argu_contra}
\tilde{v}(\tau)\leq v(\tau)\quad \mbox{for every}\quad \tau\geq 0.
\end{equation}
From inequality (\ref{ergo_eq1}) of Proposition \ref{ergo}, we have:
\begin{equation}\label{argu_contra1}
|v(\tau)-\l\tau|\leq \xi\quad \mbox{and}\quad
|\tilde{v}(\tau)-\tilde{\l}\tau|\leq \xi\quad \mbox{for all}\quad \tau\geq 0.
\end{equation}
We then easily conclude that $\tilde{\l}\leq \l$ as a consequence of
(\ref{argu_contra}) and (\ref{argu_contra1}).\\

\noindent {\bf Step 2: Continuity of $\o{f}$}.\\

\noindent We refer the reader to Proposition \ref{lwt3b} which implies
in particular the continuity of $\o{f}$. The main idea of the proof is
to apply a perturbation argument using the inequality
$|\l^{\pm}(T)-\l|\leq \frac{\xi}{T}$. $\hfill{\Box}$\\

\noindent \textbf{Proof of Proposition \ref{ganaza}.}\\

\noindent \textit{Global existence.} This is a direct consequence of the
Cauchy-P\'{e}ano theorem, using in particular the continuity of $f$ (see
\ref{Ymala3}).\\

\noindent \textit{Uniqueness.} Assume that there exists $u^{1}\in
C^{1}([0,\infty);\R)$ another solution of (\ref{homog_eq}). Define
$k(t)=|u^{0}(t)-u^{1}(t)|$, we compute (with the sign function
$\sgn(x)=x/|x|$ if $x\neq 0$):
\begin{eqnarray*}
k_{t}(t)&=&(u^{0}_{t}(t)-u^{1}_{t}(t))\,\sgn(u^{0}(t)-u^{1}(t))\\
&=& (\o{f}(u^{0}(t),t)-\o{f}(u^{1}(t),t))\,\sgn(u^{0}(t)-u^{1}(t))\leq 0,
\end{eqnarray*}
where for the last line we have used the monotonicity of $\o{f}$
(see (\ref{Ymala3})). This immediately implies that
$u^{0}=u^{1}$. $\hfill{\Box}$

%
%
%
%
%

\section{A stability result for the effective slope $\o{f}$}\label{sec3}
In this section, we will show a stability result for the term $\l$ given by
Proposition \ref{ergo} under a perturbation of (\ref{ergo_eq}) of
the form:
\begin{equation}\label{pert_ergo_eq}
\left\{
\begin{aligned}
&v_{t}=g_{\g}(v,t)=g(v,t)\pm\g,\quad t>0\\
&v(0)=v_{0},
\end{aligned}
\right.
\end{equation}
where $\g>0$ is a small real number. More precisely, we have the
following proposition:
\begin{proposition}\textbf{(Stability result).}\label{prop_pert}
Take $0< \g<1$. Let $\l_{\g}$ be the effective slope given by Proposition
\ref{ergo}, which is associated to equation (\ref{pert_ergo_eq}), and let
$\l_{0}$ be the one corresponds to $\g=0$ in (\ref{pert_ergo_eq}).
Then we have the following estimate:
\begin{equation}\label{lambda_ext}
|\l_{\g}-\l_{0}|\leq \frac{\bar{\xi}}{|\log\g|}\quad \mbox{with}\quad
\bar{\xi}=(3+2\xi)(1+2L), 
\end{equation}
and $L=\left\|\frac{\partial g}{\partial v}\right\|_{L^{\infty}(\R^{2})}$.
\end{proposition}
\textbf{Proof.} We assume, for the
sake of simplicity, that $g_{\g}=g+\g$. In this case, it is easy to
check that $\l_{\g}\geq \l_{0}$. The other case with $g_{\g}=g-\g$ is treated
similarly. We first transform our ODE problem into a PDE one by
setting $v^{\g}$ as the solution of the following equation:
\begin{equation}\label{ode-pde}
\left\{
\begin{aligned}
&v^{\g}_{t}(t,x)=g(v^{\g}(t,x),t)+\g,&\mbox{in}\quad(0,\infty)\times \R\\
&v^{\g}(0,x)=x,&x\in \R.
\end{aligned}
\right.
\end{equation}
The proof is divided into three steps.\\

\noindent \textbf{Step 1: A control on $v^{\g}_{x}$}.\\

\noindent Using comparison principle arguments for (\ref{ode-pde}), it is
easily checked that $v^{\g}(t,.)$ is a non-decreasing function
satisfying $v^{\g}(t,x+1)=v^{\g}(t,x)+1$. We want to control
$v^{\g}_{x}(t,.)$ for any $t$. For this reason, we proceed as
follows. Define for $z\geq 0$:
$$
\eta(t,x)=v^{\g}(t,x+z)-ze^{Lt},\quad t>0,\;x\in \R.
$$
We compute
\begin{eqnarray*}
\eta_{t}(t,x)&=&v^{\g}_{t}(t,x+z)-zLe^{Lt}\\
&=& g(\eta(t,x)+ze^{Lt},t)-zLe^{Lt}+\g\\
&\leq& g(\eta(t,x),t)+\g,
\end{eqnarray*}
which proves that $\eta$ is a sub-solution of (\ref{ode-pde}) with
$\eta(0,x)=v^{\g}(0,x+z)-z=v^{\g}(0,x)$, and therefore, by
the comparison principle, we obtain
$$\eta(t,x)=v^{\g}(t,x+z)-ze^{Lt}\leq v^{\g}(t,x)$$
hence for any $t\geq 0$, we have
$0\leq v^{\g}(t,x+z)-v^{\g}(t,x)\leq ze^{Lt}$,
then $v^{\g}(t,x)$ is Lipschitz continuous in the variable $x$, satisfying:
\begin{equation}\label{h1}
0\leq v^{\g}_{x}(t,x)\leq e^{Lt}\quad \mbox{ for } t\geq 0 \mbox{ and a.e. } x
\in \R.
\end{equation}
In a similar way, we can obtain a positive bound from below on
$v^{\g}_{x}$, and finally get
\begin{equation}\label{formal}
e^{-Lt}\leq v^{\g}_{x}(t,x)\leq e^{Lt}.
\end{equation}

\noindent {\bf Step 2: An upper bound of $v^{\g}$}.\\

\noindent We seek to find an upper bound of $v^{\g}$ by constructing
an explicit super-solution of (\ref{ode-pde}) with suitable initial data,
and comparing it with $v^{\g}$. For this purpose, let
\begin{equation}\label{w}
w(t,x)=v^{0}(t,x+c_{1}\g t),\quad (t,x)\in (0,\infty)\times \R,
\end{equation}
where $v_{0}$ is equal to $v^{\g}$ for $\g=0$, and $c_{1}$ is
positive constant to be precised later. We calculate: 
\begin{eqnarray*}
w_{t}(t,x)&=&v^{0}_{t}(t,x+c_{1}\g t)+c_{1}\g v^{0}_{x}(t,x+c_{1}\g t)\\
&=& g(w(t,x),t)+c_{1}\g v^{0}_{x}(t,x+c_{1}\g t),
\end{eqnarray*}
where from (\ref{formal}), we deduce that
\begin{equation}\label{est_w}
 w_{t}(t,x) \geq g(w(t,x),t)+c_{1}\g e^{-Lt}.
\end{equation}
Take $c_{1}=e^{LT}$ for some fixed $T>0$. Then using (\ref{est_w}), we
get $w_{t}(t,x) \geq g(w(t,x),t)+\g$ for any $t\in [0,T]$.
Hence $w$ is a super-solution of (\ref{ode-pde}) over $[0,T]$ whose
initial condition $w(0,x)=v^{\g}(0,x)$, which finally gives:
\begin{equation}\label{step2}
w(t,x)\geq v^{\g}(t,x)\quad \forall t\in [0,T],\;x\in \R.
\end{equation}

\noindent \textbf{Step 3: Conclusion}.\\

\noindent We will now show the error estimate (\ref{lambda_ext}). To
this end, we will estimate both sides of inequality (\ref{step2})
involving $\l_{0}$ and $\l_{\g}$. Firstly, using (\ref{ergo_eq1}) and
(\ref{formal}), we compute: 
\begin{eqnarray*}
|v^{0}(t,x+e^{LT}\g t)-v^{0}(0,x)|&\leq& |v^{0}(t,x+e^{LT}\g
t)-v^{0}(t,x)|+|v^{0}(t,x)-v^{0}(0,x)|\\
&\leq& e^{Lt}e^{LT}\g t+\l_{0}t+\xi.
\end{eqnarray*}
We take this inequality for $t=T$ and $x=0$, we get
\begin{equation}\label{s1}
w(T,0)=v^{0}(T, e^{LT}\g T)\leq \g Te^{2LT} +\l_{0}T+\xi.
\end{equation}
Secondly, using similar arguments, and the fact that $\g<1$, we
obtain $|v^{\g}(T,0)-v^{\g}(0,0)-\l_{\g}T|\leq 2+\xi$, hence
\begin{equation}\label{s2}
v^{\g}(T,0)\geq \l_{\g}T-(2+\xi).
\end{equation}
Combining (\ref{step2}), (\ref{s1}) and (\ref{s2}), it follows that
\begin{equation}\label{s3}
(\l_{\g}-\l_{0})T\leq \g Te^{2LT} + 2(1+\xi)
\end{equation}
Using (\ref{s3}), we deduce that:
\begin{equation}\label{s4}
|\l_{\g}-\l_{0}|\leq \g e^{2LT}+\frac{2(1+\xi)}{T}.
\end{equation}
Since the variable $T$ was arbitrary chosen, let $T$ satisfies
$\g T e^{2LT}=1$ and therefore $T\geq \frac{|\log \g|}{1+2L}$.
From (\ref{s4}), the result directly follows.
$\hfill{\Box}$\\

\noindent An immediate consequence of Proposition \ref{prop_pert} is the
following:
\begin{proposition}\textbf{(Modulus of continuity of $\o{f}$).}\label{lwt3b}
The function $\o{f}(u,t)$ given by (\ref{Ymala2}) satisfies for any
$(u,t)\in \R^{2}$, and for all $|v|+|s|< \frac{1}{\a}$:
\begin{equation}\label{lwt3b1}
|\o{f}(u+v,t+s)-\o{f}(u,t)|\leq \frac{\bar{\xi}}{|\log \a(|v|+|s|)|},
\end{equation}
where $\a$ is given in assumption (A1).
\end{proposition}
Remark that estimate (\ref{lwt3b1}) is optimal in view of Example \ref{example2}.

%
%
%
%
%
%

\section{Basic error estimate}\label{sec4} 
We start this section by considering a discrete scheme associated to the ODE
(\ref{homog_eq}). Namely, for a given $v^{0}$ (which may be chosen equal
to $u_{0}$ or may be different), and for a time step $\Delta t>0$, we
define the sequence $(v^{k})_{k\in \N}$ as follows: 
\begin{equation}\label{dis_sch}
v^{k+1}=v^{k}+\l_{k}\D t,\quad \l_{k}=\o{f}(v^{k},k\D t),\quad k\in \N.
\end{equation}
In this section we give a local
error estimate between the solution $u^{\eps}$ of (\ref{main_eq}) and
the sequence $v^{k}$. To be more precise, we will show the following proposition:
\begin{proposition}\textbf{(Basic error estimate).}\label{bef_proof}
  Under assumptions (A1), (A2) and
(A3) on the function $f$, let $\D t>0$ be small enough (depending only
on $\a$ and $\beta$), and take
$$e_{i}=|u^{\eps}(i \D t)-v^{i}|,\quad i=0,1.$$
Then we have:
\begin{equation}\label{sh7mra}
e_{1}\leq e_{0}+c\eps + \frac{c \D t}{|\log c \D t|},
\end{equation}
where $c=c(\a,\beta)>0$ is a positive constant.
\end{proposition}
The proof of the above proposition will be presented later in this
section. In what follows in this section and in Section \ref{sec5}, we
will assume that the arguments of the various logarithms are all less
than $1$.
\begin{lemma}\textbf{(Refined basic error estimate).}\label{lem1} 
Under the same hypothesis of Proposition 
\ref{bef_proof}, let
$$e_{i}^{+}=\max(0,u^{\eps}(i \D t)-v^{i}),\quad
e_{i}^{-}=\min(0,u^{\eps}(i \D t)-v^{i}),\quad i=0,1$$
and let $\displaystyle d^{+}_{1}= \sup_{t\in [0,\D t]}(\max (0,
u^{\eps}(t)-v^{0}))$, 
$\displaystyle d_{1}^{-}=\inf_{t\in [0,\D t]}(\min (0, u^{\eps}(t)-v^{0}))$.
Then we have:
\begin{equation}\label{sir_d1}
e^{+}_{1}\leq e_{0}^{+} + (2+\xi)\eps + \frac{\bar{\xi} \D t}{|\log \a
(|d^{-}_{1}|+ \D t)|}
\end{equation}
and
\begin{equation}\label{sir_d2}
e^{-}_{1}\geq e^{-}_{0} - (2+\xi)\eps - \frac{\bar{\xi} \D t}{|\log \a
(d^{+}_{1}+\D t)|},
\end{equation}
with $\xi=1+2\beta$ and $\beta$ defined in (A1).
\end{lemma}
\textbf{Proof.} In order to get estimates (\ref{sir_d1}),
(\ref{sir_d2}), the main idea is to freeze the last two arguments of
$f$, and to use some comparison arguments. We start by estimating the
term $f\left(\frac{u^{\eps}}{\eps},\frac{t}{\eps},u^{\eps},t\right)$ from
below. Since $u^{\eps}(t)\leq v^{0}+d_{1}^{+}$ for $t\in [0,\D t]$, we
deduce, using (A3), that 
\begin{equation*}
f\left(\frac{u^{\eps}}{\eps},\frac{t}{\eps},u^{\eps},t\right)-
f\left(\frac{u^{\eps}}{\eps},\frac{t}{\eps},v^{0},0\right)
\geq
f\left(\frac{u^{\eps}}{\eps},\frac{t}{\eps},v^{0}+d^{+}_{1},t\right)
-f\left(\frac{u^{\eps}}{\eps},\frac{t}{\eps},v^{0},0\right),
\end{equation*}
and hence, from (A1), we get
\begin{equation}\label{sir_d3}
f\left(\frac{u^{\eps}}{\eps},\frac{t}{\eps},u^{\eps},t\right)\geq
f\left(\frac{u^{\eps}}{\eps},\frac{t}{\eps},v^{0},0\right)-\a(d^{+}_{1}+
\D t).
\end{equation}
Using similar arguments, we can also show
\begin{equation}\label{sir_d4}
f\left(\frac{u^{\eps}}{\eps},\frac{t}{\eps},u^{\eps},t\right)\leq
f\left(\frac{u^{\eps}}{\eps},\frac{t}{\eps},v^{0},0\right)+\a(|d^{-}_{1}|+
\D t).
\end{equation}
We know from the definition of $e^{+}_{0}$ and $e_{0}^{-}$ that:
\begin{equation}\label{sir_d5}
v^{0}+e^{-}_{0} \leq u^{\eps}(0) \leq v^{0}+ e^{+}_{0}.
\end{equation}
Let $\overline{w}^{\eps}$ and $\underline{w}^{\eps}$ be the solutions of
the following 
ODEs:
\begin{equation}\label{Jin1}
\left\{
\begin{aligned}
&
\overline{w}^{\eps}_{t}=f\left(\frac{\overline{w}^{\eps}}{\eps},\frac{t}{\eps},v^{0},0\right)-\a(d^{+}_{1}+
\D t),\quad t>0\\
&\overline{w}^{\eps}(0)=v^{0}+e^{-}_{0}
\end{aligned}
\right.
\end{equation}
and
\begin{equation}\label{Jin2}
\left\{
\begin{aligned}
&
\underline{w}^{\eps}_{t}=f\left(\frac{\underline{w}^{\eps}}{\eps},\frac{t}{\eps},v^{0},0\right)+\a(|d^{-}_{1}|+ 
\D t),\quad t>0\\
&\underline{w}^{\eps}(0)=v^{0}+e^{+}_{0},
\end{aligned}
\right.
\end{equation}
respectively. From (\ref{sir_d3}), (\ref{sir_d4}) and (\ref{sir_d5}), we deduce,
using the comparison principle, that:
\begin{equation}\label{sir_d6}
\overline{w}^{\eps}\leq u^{\eps}\leq \underline{w}^{\eps}\quad
\mbox{on}\quad [0,\D t].
\end{equation}
Applying Proposition \ref{ergo} to the functions
$\overline{w}^{\eps}$ and $\underline{w}^{\eps}$, we
know that there exists two real numbers $\overline{\l}_{0}$ and
$\underline{\l}_{0}$ such that for $t\in [0,\D t]$ (assuming $\a
(d_{1}^{+}+\D t)\leq 1$ and $\a(|d_{1}^{-}|+\D t)\leq 1$):
\begin{equation}\label{sir_d7}
|\overline{w}^{\eps}(t)-\overline{w}^{\eps}(0)-\overline{\l}_{0}
t|\leq (2+\xi)\eps\quad \mbox{and}\quad
|\underline{w}^{\eps}(t)-\underline{w}^{\eps}(0)-\underline{\l}_{0}
t|\leq (2+\xi)\eps.
\end{equation}
Inequalities (\ref{sir_d6}) and (\ref{sir_d7}) give:
\begin{equation}\label{sir_d8}
v^{0}+e^{-}_{0}+ \overline{\l}_{0} t - (2+\xi)\eps \leq
u^{\eps}(t)\leq v^{0}+e^{+}_{0}+\underline{\l}_{0}t+
(2+\xi)\eps,\quad  t\in [0,\D t].
\end{equation}
Using Proposition \ref{prop_pert}, we obtain, for
$\l_{0}=\o{f}(v^{0},0)$, that
\begin{equation}\label{Jom1}
0 \geq \overline{\l}_{0}-\l_{0}\geq -\frac{\bar{\xi}}{|\log \a (d_{1}^{+} + \D
  t)|}
\end{equation}
and
\begin{equation}\label{Jom2}
0\leq \underline{\l}_{0}-\l_{0}\leq \frac{\bar{\xi}}{|\log \a (|d^{-}_{1}|+\D
  t)|}.
\end{equation}
The above two inequalities, together with (\ref{sir_d5}) and
(\ref{sir_d8}) give the result. $\hfill{\Box}$\\

Two immediate corollaries of the above lemma are the following:
\begin{corollary}\textbf{(Refined estimates involving
    $e^{\pm}_{i}$).}\label{corr1} 
Under the same hypothesis of Lemma \ref{lem1}, we have:
\begin{equation}\label{eqn_cor1}
e^{+}_{1}\leq e_{0}^{+} + (2+\xi)\eps + \frac{\bar{\xi} \D t}{|\log \a_{1}
(|e^{-}_{0}|+ \D t)|}
\end{equation}
and
\begin{equation}\label{eqn_cor2}
e^{-}_{1}\geq e^{-}_{0} - (2+\xi)\eps - \frac{\bar{\xi} \D t}{|\log \a_{1}
(e_{0}^{+}+\D t)|},
\end{equation}
where $\a_{1}=\a_{1}(\a,\beta)>0$ is a positive constant.
\end{corollary}
\textbf{Proof.} We have $d^{+}_{1}\leq e^{+}_{0}+\beta\D t$ and $|d^{-}_{1}|\leq
|e^{-}_{0}|+\beta\D t$ (recall that $\beta=\|f\|_{L^{\infty}(\R^{4})}$).
Therefore, the result can be deduced from (\ref{sir_d1}) and
(\ref{sir_d2}). $\hfill{\Box}$
\begin{corollary}\textbf{(Refined estimates with continuous time).}\label{corr2}
Under the same hypothesis of Proposition \ref{bef_proof}, for every
$t\in [0, \D t]$, define the continuous function $\bar{e}$ by:
\begin{equation}\label{ebar}
\bar{e}(t)=u^{\eps}(t)-(v^{0}+\l_{0} t).
\end{equation}
Also define $e^{+}(t)=\max (0, \bar{e}(t))$ and $e^{-}(t)=\min (0,
\bar{e}(t))$. Then we have for $0\leq t_{1}<t_{2}\leq \D t$ and $c=c(\a,\beta)>0$:
\begin{equation}\label{79bas}
e^{+}(t_{2})\leq e^{+}(t_{1}) + c\eps + \frac{c (t_{2}-t_{1})}{|\log
c (|e^{-}(t_{1})|+ (t_{2}-t_{1}))|}
\end{equation}
and
\begin{equation}\label{79basm}
e^{-}(t_{2})\geq e^{-}(t_{1}) - c\eps - \frac{c (t_{2}-t_{1})}{|\log
c (e^{+}(t_{1})+(t_{2}-t_{1}))|}.
\end{equation}
\end{corollary}
\textbf{Proof.} We apply the same proof (word by word) of Lemma
\ref{lem1} and Corollary \ref{corr1}, with the origin $0$ shifted to
$t_{1}$ and the time step $\D t$ replaced by $t_{2}-t_1$. $\hfill{\Box}$\\

Now we are ready to prove Proposition \ref{bef_proof}.\\

\noindent \textbf{Proof of Proposition \ref{bef_proof}.} From the
definition (\ref{ebar}) of $\bar{e}$, we know that
$$e_{0}= |\bar{e}(0)|\quad \mbox{and}\quad e_{1}=|\bar{e}(\D t)|.$$
The case where $e_{1}=\bar{e}(\D t)=0$ is obvious. Four cases could be considered.\\

\noindent  \textit{\textbf{Case 1.}} \textit{($\bar{e}(0)\geq 0$ and
$\bar{e}(\D t)>0$)}. In this case we have $e_{0}=e^{+}_{0}$,
$e_{1}=e^{+}_{1}$ and $e^{-}_{0}=0$. Therefore (\ref{sh7mra}) is an
immediate
consequence of (\ref{eqn_cor1}).\\

\noindent  \textit{\textbf{Case 2.}} \textit{($\bar{e}(0)\leq 0$ and
$\bar{e}(\D t) < 0$)}. Similar to Case 1.\\

\noindent  \textit{\textbf{Case 3.}} \textit{($\bar{e}(0)< 0$ and
$\bar{e}(\D t) > 0$)}. In this case $e_{1}=e^{+}(\D t)$. Let
 the time $t_{-+}$ be defined as follows:
$$t_{-+}=\max \{t\in (0,\D t);\quad \bar{e}(t)=0 \}.$$
Using inequality (\ref{79bas}) with
$t_{1}=t_{-+}$, $t_{2}=\D t$ and $\underline{\D} t=\D t-t_{-+}\leq
\D t$, we get:
$$e_{1}\leq c\eps + \frac{c \underline{\D}t}{|\log c\underline{\D}t|}\leq
e_{0}+c\eps + \frac{c \D t}{|\log c\D t|}$$ where we have used
the fact that $e^{+}(t_{1})=e^{-}(t_{1})=0$, and hence
(\ref{sh7mra}) follows.\\

\noindent  \textit{\textbf{Case 4.}} \textit{($\bar{e}(0)> 0$ and
$\bar{e}(\D t) < 0$)}. Similar to Case 3. $\hfill{\Box}$

%
%
%
%
%
%

\section{Estimate of the rate of convergence}\label{sec5}
This section is entirely devoted to the proof of Theorem
\ref{theo2}. Let $T>0$ and let $\D t>0$ be such that
\begin{equation}\label{Reb1}
n \D t=T, \quad n\in \N,
\end{equation}
where $n$ to be chosen large enough (the choice of $\D t$ will be
given later). In order to estimate
$\|u^{\eps}-u^{0}\|_{L^{\infty}(0,T)}$, we add and subtract
$v$, the continuous piecewise linear function passing through the points $(k\D t,
v^{k})$, $k=0\cdots n$. In other words
$$
v(t)=v^{k}+(t-k\D t)\l_{k},\quad k\D t\leq t\leq (k+1)\D t,\quad
k=0\cdots n-1,
$$
where $v^{k}$ and $\l_{k}$ are defined in (\ref{dis_sch}). 
We start by stating the following corollary that generalizes Proposition
\ref{bef_proof}.
\begin{corollary}\textbf{(Basic error estimate in continuous
    time).}\label{gene_bef} 
Under assumptions (A1), (A2) and (A3), let $T>0$, and $\D t>0$ given
by (\ref{Reb1}). Define the function $e(t)$ by:
\begin{equation}\label{Reb2}
e(t)=|u^{\eps}(t)-v(t)|,\quad t\in [0,T].
\end{equation}
For $k=0\cdots n-1$, call $e_{k}=e(k\D t)$. Then for $k\D t\leq
t\leq (k+1)\D t$, we have:
\begin{equation}\label{Reb3}
e(t)\leq e_{k}+c\eps +\frac{c\D t}{|\log c \D t|},
\end{equation}
where $c=c(\a,\beta)>0$ is a positive constant (the same given by
Proposition \ref{bef_proof}).
\end{corollary}
\noindent \textbf{Proof.} We simply apply Proposition \ref{bef_proof}
with the origin $0$ shifted to $k\D t$ and with the time step $\D t$
replaced by $\delta t = t-k\Delta t \leq \Delta t$. $\hfill{\Box}$\\

At this stage, we can show an inequality similar to (\ref{Reb3}) with
$e(t)$ and $e_{k}$ replaced respectively by the functions
\begin{equation}\label{Reb89}
e^{0}(t)=|u^{0}(t)-v(t)|\quad \mbox{and}\quad e^{0}_{k}=e^{0}(k\D t),
\end{equation}
where $u^{0}$ is the solution of
(\ref{homog_eq}). Indeed, we have the following proposition:
\begin{proposition}\textbf{(Basic error estimate for the homogenized
    ODE).}\label{mitl_f} 
Let $\o{f}$ be the function given by (\ref{Ymala2}), and enjoying the
properties given by Propositions \ref{AM7Z} and \ref{lwt3b}. Let
$T>0$, and $\D t>0$ given by (\ref{Reb1}). Then for $k\D t\leq
t\leq (k+1)\D t$, $k=0\cdots n-1$, we have:
\begin{equation}\label{Reb3bis}
e^{0}(t)\leq e^{0}_{k}+\frac{c\D t}{|\log c\D t|}
\end{equation}
where $e^{0}$, $e^{0}_{k}$ are given by (\ref{Reb89}), and
$c=c(\a,\beta)>0$ is a positive constant. 
\end{proposition}
\textbf{Sketch of the proof.} Although the proof is an adaptation of
the proof of inequality (\ref{Reb3}) (it
suffices to deal with $\o{f}(u,t)$ instead of $f(v,\t,u,t)$, and to
take $\eps=0$), we will indicate the main points where it slightly
differs. The crucial idea is that, on the one hand, the monotonicity of
the function $f(v,\t,u,t)$ with respect to the variable $u$ (see
assumption (A3)) is replaced by the monotonicity of $\o{f}(u,t)$
with respect to $u$ (see (\ref{Ymala3})). On the
other hand, the fact that $f$ is Lipschitz continuous (see
assumption (A1)) is replaced by the fact that $\o{f}$ satisfies a
modulus of continuity (see Proposition \ref{lwt3b}), and the fact
that $\o{f}(u,t)=\l(u,t)$.

Let us go into the details. In fact, the proof of Lemma \ref{lem1} can be
adapted where inequalities (\ref{sir_d3}) and
(\ref{sir_d4}) are 
replaced by
\begin{equation}\label{Gata1}
\o{f}(u^{0},t)\geq \o{f}(v^{0},0)-\frac{\bar{\xi}}{|\log
  \a(d^{0,+}_{1}+\D t)|},\quad \mbox{with}\quad d^{0,+}_{1}=\sup_{t\in
    [0,\D t]}(\max(0,u^{0}(t)-v^{0} ))
\end{equation}
and
\begin{equation}\label{Gata2}
\o{f}(u^{0},t)\leq \o{f}(v^{0},0)+\frac{\bar{\xi}}{|\log
  \a(|d^{0,-}_{1}|+\D t)|},\quad \mbox{with}\quad d^{0,-}_{1}=\inf_{t\in
    [0,\D t]}(\min(0,  u^{0}(t)-v^{0}))
\end{equation}
respectively. Here we have used the monotonicity of $\o{f}$, and inequality (\ref{lwt3b1}) given by
Proposition \ref{lwt3b}. Having (\ref{Gata1}) and (\ref{Gata2}) in
hands, the sub- and super-solution $\overline{w}^{\eps}$,
$\underline{w}^{\eps}$ defined by (\ref{Jin1}) and (\ref{Jin2}) are
replaced by $\overline{w}^{0}$, $\underline{w}^{0}$ solutions of
$$
\left\{
\begin{aligned}
& \overline{w}^{0}_{t}=\o{f}(v^{0},0)-\frac{\bar{\xi}}{|\log
  \a(d^{0,+}_{1}+\D t)|}\\
&\overline{w}^{0}(0)=v^{0}+e^{0,-}\quad\mbox{with}\quad e^{0,-}=\min(0,u^{0}(0)-v^{0})
\end{aligned}
\right.
$$
and
$$
\left\{
\begin{aligned}
& \underline{w}^{0}_{t}=\o{f}(v^{0},0)+\frac{\bar{\xi}}{|\log
  \a(|d^{0,-}_{1}|+\D t)|}\\
&\underline{w}^{0}(0)=v^{0}+e^{0,+}\quad\mbox{with}\quad e^{0,+}=\max(0,u^{0}(0)-v^{0})
\end{aligned}
\right.
$$
respectively, and we have
\begin{equation}\label{Gata3-}
\overline{w}^{0}(t)-\overline{w}^{0}(0)-\overline{\l}^{0}_{0}t=0\quad
\mbox{and} \quad
\underline{w}^{0}(t)-\underline{w}^{0}(0)-\underline{\l}^{0}_{0}t=0
\end{equation}
with
\begin{equation}\label{Gata3}
\overline{\l}^{0}_{0}=\overline{w}^{0}_{t}=\l_{0}-\frac{\bar{\xi}}{|\log
  \a(d^{0,+}_{1}+\D t)|}\quad \mbox{and} \quad
\underline{\l}^{0}_{0}=\underline{w}^{0}_{t}=\l_{0}+\frac{\bar{\xi}}{|\log
  \a(|d^{0,-}_{1}|+\D t)|},
\end{equation}
where we recall the reader that $\l_{0}=\o{f}(v^{0},0)$. 
Remark that (\ref{Gata3-}) replaces (\ref{sir_d7}), while (\ref{Gata3})
gives (as a replacement of (\ref{Jom1}) and (\ref{Jom2})):
$$0\geq \overline{\l}^{0}_{0}-\l_{0}=-\frac{\bar{\xi}}{|\log
  \a(d^{0,+}_{1}+\D t)|}\quad \mbox{and}\quad 0\leq
\underline{\l}^{0}_{0}-\l_{0}=\frac{\bar{\xi}}{|\log
  \a(|d^{0,-}_{1}|+\D t)|}.$$
At this point, the rest of the proof proceeds in a very similar way as
the proof of (\ref{Reb3}) with $\eps=0$, and a possible changing of
the constants but always depending on $\a$ and $\beta$. $\hfill{\Box}$\\

Now we are ready to present the proof of our main result.\\

\noindent \textbf{Proof of Theorem \ref{theo2}.} We first note that
the constant $c>0$ may certainly differ from line to line in the
proof. We decompose the quantity
$|u^{\eps}(t)-u^{0}(t)|$ into two pieces:
\begin{equation}\label{Reb7}
|u^{\eps}(t)-u^{0}(t)|\leq \overbrace{|u^{\eps}(t)-v(t)|}^{e(t)}+
\overbrace{|u^{0}(t)-v(t)|}^{e^{0}(t)}.
\end{equation} In order to
estimate $e(t)$, we iterate inequality (\ref{Reb3}) and we finally
obtain for $u^{\eps}(0)=u_{0}=v^{0}$, and $t\in [0,T]$ with $T=n\D t$:
\begin{equation*}
e(t) \leq e_{0}+ cn\eps + \frac{c n \D t}{|\log c\D t|}\leq \frac{c
  \eps T }{\D t}+ \frac{c T}{|\log c \D t|}.
\end{equation*}
The above inequality gives:
\begin{equation}\label{Reb4}
e(t)\leq cT\left(\frac{\eps}{\D t}+\frac{1}{|\log \D t|} \right),
\end{equation}
and, from inequality (\ref{Reb3bis}) of Proposition \ref{mitl_f}, we can
show in the same way as above that we also have:
\begin{equation}\label{Reb5}
e^{0}(t)\leq \frac{cT}{|\log \D t|}.
\end{equation}
Choosing particularly $\D t = C \eps |\log \eps|$, we deduce that
$\eps \ll \D t = C \eps |\log \eps| \ll T $, and (from (\ref{Reb4}),
(\ref{Reb5})) that:
\begin{equation}\label{Reb6}
e(t)\leq \frac{c T}{|\log \eps|}\quad \mbox{and}\quad e^{0}(t)\leq
\frac{c T}{|\log \eps|},
\end{equation}
with $c$ in (\ref{Reb6}) depending on the choice of $C>0$. finally,
inequality (\ref{error}) could now be 
easily deduced from (\ref{Reb7}) and (\ref{Reb6}). $\hfill{\Box}$

%
%
%
%

\section{Application: error estimate for linear transport
  equations}\label{sec6}
In this section, as an application of our previous results on ODEs, we give
the proof of some error estimates for the homogenization of linear
transport equations. Namely, we prove Theorems \ref{theo3} and
\ref{theo4}. We start with Theorem \ref{theo3} keeping the same
notations of Subsection \ref{sub1.2}.\\

\noindent\textbf{Proof of Theorem \ref{theo3}.} The proof is divided
into four steps. The first three steps are devoted to the definition of
the limit solution $V^{0}$, and to prove that it is a viscosity
solution. The proof of the error estimate is done in the last step.\\

\noindent {\bf Step 1: Definition of $V^{0}$}.\\

\noindent Because the homogenized vector field $\o{a}$ is not Lipschitz,
we define our solution $V^{0}$ to (\ref{transport_eqn0}) in an indirect
way using the characteristics. Precisely, for $(t,x)\in (0,\infty)\times
\R^{2}$, $x=(x_{1},x_{2})$, we 
define $V^{0}(t,x)$ as follows: 
\begin{equation}\label{eq1_sec6}
V^{0}(t,x)=V_{0}(X^{0}(0;t,x))
\end{equation}  
where the curve $X^{0}(\t;t,x):\t\in \R\rightarrow X^{0}(\t;t,x)\in
\R^{2}$,
\begin{equation}\label{fsof_s}   
X^{0}(\t;t,x)=(X_{1}^{0}(\t;t,x), X_{2}^{0}(\t;t,x))
\end{equation}
is the solution of the following ODE:
\begin{equation}\label{eq2_sec6}
\left\{
\begin{aligned}
&\partial_{\t} X^{0}=\o{a} (X^{0})\\
&X^{0}(t)=x. 
\end{aligned}
\right.
\end{equation}
We will see below that this solution is unique. For the sake of
simplicity of notations, we will omit the dependence of 
$X^{0}$ on $(t,x)$ and we will simply write 
$$X^{0}(\t;t,x)=X^{0}(\t).$$
From (\ref{eq2_sec6}) and (\ref{att_sign}), we can easily check that
$X^{0}_{2}(\t)=x_{2}-t+\tau$, hence using (\ref{eq1_sec6}), we get
$$V^{0}(t,x)=V_{0}\left(X_{1}^{0}(0), x_{2}-t\right)$$ where $X_{1}^{0}$ satisfies:
\begin{equation}\label{eq3_sec6}
\left\{
\begin{aligned}
&\partial_{\t} X_{1}^{0}=-\o{f}\left(X_{1}^{0}, x_{2}-t+\tau\right)\\
&X_{1}^{0}(t)=x_{1}. 
\end{aligned}
\right.
\end{equation}
In order to show that $X_{1}^{0}(0)$ is uniquely defined, we solve
(\ref{eq3_sec6}) backwards, in other words, we let 
$$\o{X}_{1}^{0}(\t)=X_{1}^{0}(t-\t).$$
In this case:
\begin{equation}\label{7rmhl}
V^{0}(t,x)=V_{0}(\o{X}^{0}_{1}(t), x_{2}-t),
\end{equation}
where $\o{X}_{1}^{0}$ satisfies:
\begin{equation}\label{eq4_sec6}
\left\{
\begin{aligned}
&\partial_{\t} \o{X}_{1}^{0}=\o{f}\left(\o{X}_{1}^{0}, x_{2}-\tau\right)\\
&\o{X}_{1}^{0}(0)=x_{1}.
\end{aligned}
\right.
\end{equation} 
From Proposition \ref{ganaza}, the solution $\o{X}_{1}^{0}\in
C^{1}([0,\infty);\R)$ is unique and hence
$X_{1}^{0}(0)=\o{X}_{1}^{0}(t)$ is uniquely determined. Consequently the
function $V^{0}$ is well defined.\\

\noindent {\bf Step 2: $V^{0}$ is Lipschitz continuous}.\\

\noindent From Step 1, we know that 
\begin{equation}\label{Ya1}
V^{0}(t,x_{1},x_{2})=V_{0}(\o{X}_{1}^{0}(t),x_{2}-t)
\end{equation}
with $\o{X}_{1}^{0}$ given by (\ref{eq4_sec6}) also depends on
$x_{1}$ and $x_{2}$. Let the function $Y:\R^{3}\rightarrow \R$ be
defined as follows (with simplified notation showing the dependence on
the variables $(t,x_{1},x_{2})$): 
$$Y(t,x_{1},x_{2}):=\o{X}^{0}_{1}(t).$$
In order to show that $V^{0}$ is Lipschitz, it suffices (see
(\ref{Ya1})) to show that $Y$ is Lipschitz. First, it is easily seen
from (\ref{eq4_sec6}) that $Y$ is Lipschitz in time $t$. The Lipschitz continuity with
respect to the variable $x_{1}$ directly follows from the monotonicity
of $\o{f}$ (see (\ref{Ymala3})), and the comparison principle. In order
to show the Lipschitz continuity with respect to $x_{2}$, we first give a formal proof by
assuming that $\o{f}$ is smooth, and then we present the main idea that
permit to make the proof rigorous. Suppose that 
$$\o{f}\in C^{\infty}(\R^{2};\R)\quad \mbox{with}\quad |\o{f}(y,s)|\leq
\|\o{f}\|_{\infty}.$$ 
Take
$$\widehat{Y}=\partial_{x_{2}}Y\quad \mbox{and}\quad \widetilde{Y}=\partial_{t}Y,$$
then the above two functions satisfy:
\begin{equation}\label{Ya2}
\left\{
\begin{aligned}
&\partial_{t}\widehat{Y}=(\partial_{y}\o{f})\widehat{Y} + \partial_{s}\o{f}\\
&\widehat{Y}(0,x_{1},x_{2})=0, 
\end{aligned}
\right.
\end{equation} 
and
\begin{equation}\label{Ya3}
\left\{
\begin{aligned}
&\partial_{t}\widetilde{Y}=(\partial_{y}\o{f})\widetilde{Y} - \partial_{s}\o{f}\\
&\widetilde{Y}(0,x_{1},x_{2}) = \o{f}(x_{1},x_{2}), 
\end{aligned}
\right.
\end{equation} 
respectively. Let $\o{Y}=\widehat{Y}+\widetilde{Y}$, we get (from (\ref{Ya2})
and (\ref{Ya3})):
\begin{equation}\label{Ya4}
\left\{
\begin{aligned}
&\partial_{t}\o{Y}=(\partial_{y}\o{f})\o{Y}\\
&\o{Y}(0,x_{1},x_{2}) = \o{f}(x_{1},x_{2}),
\end{aligned}
\right.
\nonumber
\end{equation}
which gives, because of the monotonicity of $\o{f}$ (see (\ref{Ymala3})), 
that the function $t\rightarrow |\o{Y}(t,.,.)|$ is
non-increasing. Hence 
\begin{equation}\label{Ya4}
|\o{Y}(t,x_{1},x_{2})|\leq
|\o{f}(x_{1},x_{2})|\leq \|\o{f}\|_{\infty}.
\end{equation}   
Since $\widetilde{Y}=\partial_{t}Y$, we know from (\ref{eq4_sec6})
that $|\widetilde{Y}|\leq \|\o{f}\|_{\infty}$ where we finally obtain (see
(\ref{Ya4})):
$$|\widehat{Y}|=|\partial_{x_{2}}Y|\leq 2\|\o{f}\|_{\infty},$$
which shows that $Y$ is Lipschitz continuous in the $x_{2}$
variable. In order to make the proof rigorous, it suffices to consider a
regular approximation of the function $\o{f}$ (as for example the
convolution with a suitable mollifier sequence) and then to pass to the limit.\\  

\noindent {\bf Step 3: $V^{0}$ is a viscosity solution of
  (\ref{transport_eqn0})}.\\ 

\noindent For the definition and the study of the theory of viscosity
solutions, we refer the reader to \cite{Barles}.  Let us take $\overline{\phi},
\underline{\phi}\in C^{1}(\R^{3};\R)$ such that $V^{0}-\o{\phi}$
(resp. $V^{0}-\underline{\phi}$) has a local maximum (resp. local minimum) at
some point $(\o{t},\o{x})\in (0,\infty)\times \R^{2}$
(resp. $(\underline{t},\underline{x})\in (0,\infty)\times \R^{2}$), with
$V^{0}(\o{t},\o{x})=\o{\phi}(\o{t},\o{x})$ and
$V^{0}(\underline{t},\underline{x})=\underline{\phi}(\underline{t},\underline{x})$.
In order to show that $V^{0}$ is a viscosity solution of (\ref{transport_eqn0}), we need to
show the following two inequalities:
\begin{equation}\label{fsof_s1}
\partial_{t}\o{\phi}(\o{t},\o{x})+(\o{a}\cdot
\nabla\o{\phi})(\o{t},\o{x})\leq 0 
\end{equation}
and
\begin{equation}\label{fsof_s2}
\partial_{t}\underline{\phi}(\underline{t},\underline{x})+(\o{a}\cdot
\nabla\underline{\phi})(\underline{t},\underline{x})\geq 0. 
\end{equation}
We only show inequality (\ref{fsof_s1}). In fact, inequality
(\ref{fsof_s2}) can be proved in exactly the same way. For any $t\in
[0,\o{t}]$ define the function $\phi$ by
$$\phi:t\rightarrow \phi(t)=\o{\phi}(t,X^{0}(t;\o{t},\o{x})),$$
where $X^{0}$ is defined in (\ref{eq2_sec6}). Let us show an inequality
on $\phi$ in the interval $[\o{t}-r,\o{t}]$ for $r>0$ small enough. Remark that
$X^{0}(\o{t};\o{t},\o{x})=\o{x}$. Hence, for $t\in [\o{t}-r,\o{t}]$,
$(t,X^{0}(t;\o{t},\o{x}))$ is close to $(\o{t},\o{x})$ and therefore
(since $V^{0}-\o{\phi}$ has a local maximum at $(\o{t},\o{x})$ with
$V^{0}(\o{t},\o{x})=\o{\phi}(\o{t},\o{x})$) we get: 
\begin{eqnarray*}
& \phi(t)=\o{\phi}(t, X^{0}(t;\o{t},\o{x}))\geq
V^{0}(t,X^{0}(t;\o{t},\o{x}))=V_{0}(X^{0}(0;t,X^{0}(t;\o{t},\o{x})))\\
& = V_{0}(X^{0}(0;\o{t},\o{x}))=V^{0}(\o{t},\o{x})=\o{\phi}(\o{t},\o{x})=\o{\phi}(\o{t},  
X^{0}(\o{t};\o{t},\o{x}))=\phi(\o{t}),
\end{eqnarray*}
where the passage from the first to the second line is due to the fact
that the points $(t,X^{0}(t;\o{t},\o{x}))$ and $(\o{t},\o{x})$ are on
the same characteristics. Finally, this implies 
$$\partial_{t}\phi\big{|}_{t=\o{t}}\leq 0,$$
which directly gives (\ref{fsof_s1}).\\

\noindent {\bf Step 4: Proof of the error estimate (\ref{1.12})}.\\ 

\noindent The solution $V^{\eps}$ of (\ref{transport_eqn}) can be
written (in analogue with (\ref{7rmhl}) and (\ref{eq4_sec6})) as 
\begin{equation}\label{7rmhl_9}
V^{\eps}(t,x) = V_{0}(\o{X}_{1}^{\eps}(t), x_{2}-t)
\end{equation}
where the characteristics $\o{X}_{1}^{\eps}$ satisfies:
\begin{equation}\label{7rmhl_1}
\left\{
\begin{aligned}
&\partial_{\t}\o{X}^{\eps}_{1}=f\left(\frac{\o{X}^{\eps}_{1}}{\eps},\frac{x_{2}-\t}{\eps},
\o{X}^{\eps}_{1}, x_{2}-\t\right)\\
&\o{X}^{\eps}_{1}(0)=x_{1}.
\end{aligned}
\right.
\end{equation}
We apply Theorem \ref{theo2}, namely inequality (\ref{error}), with
$u^{\eps}$ and $u^{0}$ replaced by $\o{X}^{\eps}_{1}$ and
$\o{X}^{0}_{1}$ respectivly, we obtain:
\begin{equation}\label{7rmhl_2}
\left\|\o{X}^{\eps}_{1}-\o{X}^{0}_{1}\right\|_{L^{\infty}(0,T)}\leq \frac{cT}{|\log
  \eps|}\quad \mbox{for}\quad T\geq C\eps|\log \eps|\quad
\mbox{with}\quad C>0, \eps>0.
\end{equation}
Using (\ref{7rmhl}), (\ref{7rmhl_9}) and (\ref{7rmhl_2}), we compute for
$(t,x)\in (0,T)\times \R^{2}$:
\begin{eqnarray*}
|V^{\eps}(t,x)-V^{0}(t,x)| &=& |V_{0}(\o{X}_{1}^{\eps}(t), x_{2}-t)-
V_{0}(\o{X}_{1}^{0}(t), x_{2}-t)|\\
& \leq & Lip(V_{0})|\o{X}_{1}^{\eps}(t)-\o{X}_{1}^{0}(t)|\\
& \leq & Lip(V_{0})\frac{cT}{|\log
  \eps|},
\end{eqnarray*}
and inequality (\ref{1.12}) directly follows. $\hfill{\Box}$\\

\noindent{\textbf{Proof of Theorem \ref{theo4}.}} Application of
Proposition \ref{ergo}, and same proof as Theorem \ref{theo3}.

\section{Appendix: proof of
Examples \ref{example2} and \ref{example3}}\label{sec7}
\textbf{Proof of Example \ref{example2}.} The function $\o{f}$ can
be expressed as (see for instance \cite{Picc78}):
$$\o{f}(u,t)=\left(\int^{1}_{0}\frac{dv}{-u+|\sin 2\pi v|}\right)^{-1}.$$
Let $a=-u>0$, it is easy to check that $\displaystyle
\int^{1}_{0}\frac{dv}{a+|\sin 2\pi v|}=2\int_{|v|\leq
1/4}\frac{dv}{a+|\sin 2\pi v|}$. Take
$$ I^{a}=\int_{|v|\leq
1/4}\frac{dv}{a+|\sin 2\pi v|}\quad\mbox{and}\quad  
I^{a,R}=\int_{|v|\leq Ra}\frac{dv}{a+2\pi |v|}.$$ 
We are interested in the limit $a \rightarrow 0$ and $R\rightarrow
\infty$ with $Ra\rightarrow 0$. 
We compute
$$I^{a}-I^{a,R}= \overbrace{\int_{Ra\leq |v| \leq 1/4}\frac{dv}{a+|\sin 2\pi v|}}^{A} +
\overbrace{\int_{|v|\leq Ra} \frac{2\pi |v|-|\sin 2\pi v|}{(a+2\pi
|v|)(a+|\sin 2\pi v|)}dv}^{B},$$ where we have:
$$
\left\{
\begin{aligned}
& B \rightarrow 0\quad \mbox{as} \quad Ra\rightarrow 0\\
& A\leq \frac{1}{2}\frac{1}{\sin (2\pi R a)}\sim \frac{c}{Ra}\sim
c\sqrt{|\log a|},
\end{aligned}
\right.
$$
with $c=\frac{1}{4 \pi}$, and we have chosen 
$$R=\frac{1}{a\sqrt{|\log a|}}.$$
Now, let $\bar{v}=\frac{v}{a}$, we also
compute:
\begin{equation*}
I^{a,R} = 2 \int_{0\leq \bar{v}\leq R} \frac{d\bar{v}}{1+2\pi
\bar{v}}=\frac{1}{\pi} (\log (1+2\pi R))\sim \frac{1}{\pi} \log R \sim  \frac{1}{\pi} |\log a|.
\end{equation*}
Since $A\ll I^{a,R}$, this shows that $I^{a}\sim I^{a,R}$ and then
$$\o{f}(u,t)\sim \frac{\pi}{2|\log |u||}\quad \mbox{as} \quad
u\rightarrow 0^{-},$$
which justifies (\ref{ineq_example}).$\hfill{\Box}$\\

\noindent \textbf{Proof of Example \ref{example3}.} Since
$f(v,\t,u,t)=g(v+\t)-1$, then the function $v^{\eps}$ defined by:
$$v^{\eps}(t)=u^{\eps}(t)+t$$
satisfies
\begin{equation}\label{pr_ex3_eq1}
\left\{
\begin{aligned}
& v^{\eps}_{t}=g\left(\frac{v^{\eps}}{\eps}\right),\quad t>0,\\
& v^{\eps}(0)=0.
\end{aligned}
\right.
\end{equation}
From the particular expression (\ref{tresb}) of $g$ with
$g\left(\frac{1}{2}\right)=0$, we can check that $0\leq
\frac{v^{\eps}(t)}{\eps}\leq \frac{1}{2}$ for every $t\geq 0$, which
implies that
$$v^{\eps}\rightarrow 0\quad \mbox{in}\quad L^{\infty},$$ 
then
$$u^{\eps}\rightarrow u^{0}\mbox{ in }
L^{\infty}\quad \mbox{with} \quad u^{0}(t)=-t,$$
where
$$u_{t}^{0}=\o{f}(u^{0},t)=-1.$$
Moreover, equation (\ref{pr_ex3_eq1})
can be written:
$$v_{t}^{\eps}=\frac{1}{2}-\frac{v^{\eps}}{\eps}\quad \mbox{with}\quad
v^{\eps}(0)=0,$$
therefore (solving the above equation) we get for $t=\delta \eps |\log \eps|$
\begin{eqnarray*}
u^{\eps}(t)-u^{0}(t)=v^{\eps}(t)&=&\frac{\eps}{2}(1-e^{-\frac{t}{\eps}})\\
&=&\frac{\eps}{2}(1-e^{-\delta |\log \eps|})\\
&=&\frac{\eps}{2}-\frac{1}{2}\eps^{1+\delta}\sim \frac{\eps}{2}\sim
\frac{t}{2\delta |\log \eps|}, 
\end{eqnarray*} 
which terminates the proof of Example \ref{example3}. $\hfill{\Box}$\\

\noindent \textbf{Acknowledgements.} This work was supported by the
contract ANR MICA (2006-2009). The second author would like to thank
B. Corbin for fruitful discussions.

\bibliographystyle{siam}
\bibliography{biblio}
\end{document}